\documentclass[twoside,notitlepage,12pt]{article}

\usepackage{amssymb}
\usepackage[leqno]{amsmath}
\usepackage{amsfonts}
\usepackage{latexsym}
\usepackage{amsopn}
\usepackage{amstext}
\usepackage{amsthm}

\usepackage[all]{xy}
\newdir{ >}{{}*!/-9pt/@{>}}

\usepackage{verbatim}
\usepackage[colorlinks]{hyperref}

\providecommand{\cal}{\mathcal}

\newenvironment{pf}{\begin{proof}}{\end{proof}}

\newcommand{\Aaa}{{\cal{A}}}
\newcommand{\Bee}{{\cal{B}}}
\newcommand{\Cee}{{\cal{C}}}

\newcommand{\Ef}{{\cal{F}}}

\newcommand{\El}{{\cal{L}}}

\newcommand{\Yu}{{\cal{U}}}
\newcommand{\Vee}{{\cal{V}}}

\newcommand{\al}{\alpha}
\newcommand{\Gam}{\Gamma}

\renewcommand{\phi}{\varphi}
\renewcommand{\rho}{\varrho}

\newcommand{\rest}{\restriction}

\newcommand{\ntr}{{n\in\omega}}

\newcommand{\loe}{\leqslant}

\newcommand{\subs}{\subseteq}
\newcommand{\sups}{\supseteq}
\newcommand{\nnempty}{\ne\emptyset}

\newcommand{\ovr}{\overline}

\newcommand{\En}{\mathcal N}

\renewcommand{\S}{{\mathbb S}}

\newcommand{\cl}{\operatorname{cl}}
\newcommand{\Int}{\operatorname{int}}
\newcommand{\w}{\operatorname{w}}

\newcommand{\pr}{\operatorname{pr}}
\newcommand{\liminv}{\varprojlim}
\newcommand{\by}{/}

\newcommand{\normalsub}{\unlhd} % normal subgroup
\newcommand{\normalsup}{\unrhd} % normal supergroup
\newcommand{\diag}{\vartriangle} %\Delta

\newtheorem{tw}{Theorem}[section]
\newtheorem{wn}[tw]{Corollary}
\newtheorem{lm}[tw]{Lemma}
\newtheorem{prop}[tw]{Proposition}
\newtheorem{claim}[tw]{Claim}
\newtheorem{fact}[tw]{Fact}
\theoremstyle{definition}

\newtheorem{ex}[tw]{Example}

\newtheorem{uwga}[tw]{Remark}
\theoremstyle{remark}

\newcommand{\setof}[2]{\{#1\colon #2\}}

\newcommand{\sett}[2]{\{#1\}_{#2}}
\newcommand{\sn}[1]{\{#1\}} % singleton
\newcommand{\dn}[2]{\{#1,#2\}} % doubleton
\newcommand{\map}[3]{#1\colon #2 \to #3} % A function
\newcommand{\img}[2]{#1[#2]} % image of a set
\newcommand{\inv}[2]{{#1}^{-1}[#2]} % preimage of a set

\providecommand{\nat}{\omega}

\newcommand{\iso}{\approx}

\newcommand{\cmp}{\circ} % composition!!!

\newcommand{\invsys}[5]{\langle {#1}_{#4};{#2}_{#4}^{#5};#3 \rangle}

\newcommand{\amalgam}{\rtimes}
\newcommand{\T}{{\mathbb T}}
\newcommand{\Z}{{\mathbb Z}}

\title{A decomposition theorem for compact groups with application to supercompactness}
\author{
{\sc Wies{\l}aw Kubi\'s}
\footnote{Supported in part by the Grant IAA 100 190 901
and by the Institutional Research Plan of the Academy of Sciences of Czech Republic No. AVOZ 101 905 03.}\\
{\small Mathematical Institute}\\
{\small Czech Academy of Sciences}\\
{\small Prague, CZECH REPUBLIC}\\
{\small\texttt{wkubis@ujk.edu.pl}}
\and
{S\l awomir Turek}\\
{\small Institute of Mathematics}\\
{\small Jan Kochanowski University}\\
{\small Kielce, POLAND}\\
{\small\texttt{sturek@ujk.edu.pl}}
}

\begin{document}

\maketitle

\begin{abstract}
%The purpose of this note is to prove a decomposition theorem for connected compact groups. 
We show that every compact connected group is the limit of a continuous inverse sequence, in the category of compact groups, where each successor bonding map is either an epimorphism with finite kernel or the projection from a product by a simple compact Lie group.

As an application, we present a proof of an unpublished result of Charles Mills from 1978: every compact group is supercompact.

\vspace{3mm}


\end{abstract}

%\tableofcontents

\section{Introduction}

Several deep structure theorems for compact topological groups are known in the literature, a very good source is the book of Hofmann \& Morris \cite{HM}.
Even simple structure results allow to deduce some nontrivial topological properties of compact groups.
A good example is the theorem of Kuzminov~\cite{Kuzminov}, saying that every compact group is dyadic, i.e. a continuous image of some Cantor cube.
In particular, every disjoint family of open sets in a compact group is countable.

Trying to understand topological structure of compact groups, one immediately encounters Lie groups. In particular, it turns out that {simple compact Lie groups} form building blocks for all connected (as well as $0$-dimensional) compact groups. A \emph{simple compact Lie group} is a compact Lie group in which every proper closed normal subgroup is finite. This of course differs from the definition of a simple group in algebra. Typical examples of simple compact Lie groups are the circle group $\T$ and all finite groups.

Given a compact group $G$, by its \emph{extension} we mean a group $G_1$ together with a continuous epimorphism $\map f{G_1}G$. This notion of `extension' is better suited than the (perhaps more intuitive) usual one, where $G_1\sups G$.
In the latter definition one loses compactness when passing to the limit of infinite sequence of extensions.

There are two natural simple (from the topological point of view) extensions of compact groups: one is obtained by multiplying a given group $G$ by a simple compact Lie group
and the other is by taking an epimorphism $\map f{G_1}G$ whose kernel is finite. The second type of extension can actually lead to unexpected results, however topologically it is a local homeomorphism.
A typical (yet nontrivial) example is the map $\map f\T\T$ of the circle group, defined by $f(z)=z^2$. Here, $\ker(f)\iso\Z_2$ and the limit of an infinite sequence of these extensions gives the well known $2$-adic solenoid group.

We prove that every connected compact group can be obtained from the trivial group by using the above two types of extensions.
Namely, we prove a decomposition theorem (Theorem~\ref{ghorwgeoeqw} below) saying that every connected compact group is the limit of a continuous transfinite inverse sequence starting with the trivial group,
where each successor bonding map is an extension of one of the above mentioned types: either it has a finite kernel or it is the canonical projection from the product by a simple compact Lie group.
We also present a natural example (semidirect product of the circle $\T$ and $\Z_2$) showing that this fails for disconnected compact groups.
However, these two operations are sufficient for obtaining all Abelian as well as all $0$-dimensional compact groups.
%The key fact used in the proof is a finite decomposition of this sort, existing in every connected compact Lie group.
%It is worth mentioning that our decomposition result does not trivialize for metrizable groups.

There is one, perhaps not so well known, property of compacta, called \emph{supercompactness}.
By definition, a topological space is \emph{supercompact} if it has a subbase $\Bee$ for its closed sets, such that every linked family $\Ef\subs \Bee$ has nonempty intersection.
A family $\Ef$ is \emph{linked} if $A\cap B\nnempty$ for every $A,B\in\Ef$.
By Alexander's subbase lemma, every supercompact space is compact. A nontrivial result of Strok \& Szyma\'nski \cite{StSz} says that every compact metric space is supercompact. By the result of van Douwen \& van Mill \cite{vDvM}, every infinite supercompact space contains nontrivial convergent sequences, therefore there exist compact spaces which are not supercompact.
For more information about supercompactness we refer to van Mill's book \cite{vanMill}.

An unpublished result of Charles Mills (existing only as a seminar report \cite{Mills}) says that all compact groups are supercompact.

The final part of our note is devoted to the proof of Mills' result, applying our decomposition theorem.

\subsection{Compact Lie groups}
Recall that a \emph{compact Lie group} is a compact group $G$ which has a neighborhood of identity $U$ that does not contain nontrivial subgroups of $G$. It is a highly nontrivial fact that such a group has a manifold structure, etc. We use the above definition only for the sake of simplicity, to avoid discussing the notions and results from Lie algebra.

A \emph{simple compact Lie group} is a compact Lie group $G$ such that all its proper normal closed subgroups are finite. Note that a simple compact Lie group is either connected or finite.

We shall use the following important facts about compact groups.

\begin{fact}\label{f0}
Every closed subgroup and every continuous epimorphic image of a compact Lie group is a compact Lie group.
\end{fact}

The first part of the above fact follows immediately from our definition, while the second part is highly nontrivial, see \cite[Theorem 6.7]{HM}.

\begin{fact}[{\cite[Corollary~2.29]{HM}}]\label{f1}
Every compact group is isomorphic, in the category of topological groups, to a subgroup of a product of compact Lie groups.
\end{fact}

Thus, given a compact group $G$, we can always assume that $G\loe \prod_{i\in\Gam}H_i$, where each $H_i$ is a compact Lie group. Moreover, we may assume that $\img{\pr_i}G=H_i$ for each $i$, replacing $H_i$  by $\img{\pr_i}G$.

\begin{fact}[{\cite[Theorem~6.19]{HM}}]\label{f2}
Let $L$ be a connected compact Lie group. Then $L$ is isomorphic, in the category of topological groups, to
$$(J_0\times \dots \times J_{n-1})\by N$$
where each $J_i$, $i<n$, is a simple compact Lie group and $N$ is a finite normal subgroup of the center of $J_0\times \dots \times J_{n-1}$.
\end{fact}

\section{Main result}

Let $G$ be a compact group. A \emph{normal resolution} in $G$ is a decreasing chain $\sett{H_\al}{\al<\kappa}$
of closed normal subgroups of $G$ indexed by ordinals, such that $H_0=G$,  $H_\delta=\bigcap_{\al<\delta}H_\al$ for every limit ordinal $\delta<\kappa$ and $\bigcap_{\al<\kappa}H_\al=1$
(the trivial subgroup of $G$). Discarding the last condition, we get the notion of a \emph{normal sequence}.
Note that a normal sequence naturally induces a continuous inverse sequence of groups of the form $G\by {H_\al}$ whose limit is $G$.
On the other hand, every continuous (transfinite) inverse sequence in the category of topological groups with limit $G$ induces a normal sequence formed by kernels of the respective projections.
We shall use this fact below.

Fix a compact group $G$. We shall write $A \normalsub G$ or $G \normalsup A$ for ``$A$ is a closed normal subgroup of $G$". Writing $H_0\normalsub H_1\normalsub \dots\normalsub H$, we shall always mean that each $H_i$ is a normal subgroup of $H$ (the relation $\normalsub$ is not transitive in general).

\begin{lm}\label{wefjwei}
Let $G$ be a connected compact Lie group. Then $G$ has a normal resolution
$$G=H_0\normalsup H_1 \normalsup \dots \normalsup H_n \normalsup H_{n+1}=\{1\}$$
such that for each $i\loe n$ either $H_{i}\by {H_{i+1}}$ is finite or else there is $C\normalsub G$ such that $C\cap H_{i}=H_{i+1}$ and $G\by C$ is a simple Lie group.
\end{lm}

\begin{pf}
Let $H_0=G$. We shall use Fact~\ref{f2}, namely that $G$ is isomorphic to $J\by K$, where $J = J_0\times\dots\times J_{n-1}$, $J_0,\dots,J_{n-1}$ are simple compact Lie groups and $K$ is a finite subgroup of the center of $J$.
Let $L_i=\img{\pr_i}K$, where $\map{\pr_i}J{J_i}$ is the projection. Further, let $L=L_0\times\dots\times L_{n-1}$. Then $L$ is a finite normal subgroup of $J$, $K\subs L$, and there is a natural epimorphism $\map f{J\by K}{J\by L}$ whose kernel is isomorphic to $L\by K$, therefore finite. Identifying $J\by K$ with $G$, we set $H_{n}=\ker(f)$ and $H_{n+1}=1$.

Now observe that $J\by L = (J_0\by{L_0})\times \dots \times (J_{n-1}\by{L_{n-1}})$ and each $J_i\by{L_i}$ is a simple Lie group. 
For each $j<n$ let $$P_j= \prod_{i<j}(J_i\by{L_i})$$
and let $\map{p_j}{P_{j+1}}{P_j}$ be the canonical projection. Note that $P_0$ is the trivial group.
Finally, let
$$H_i = \ker (p_i\cmp \dots\cmp p_{n-1}\cmp f).$$
Then $H_n\by{H_{n+1}}=H_n$ is finite and for $i<n$ we have that $H_{i+1}=C_i\cap H_i$, where $C_i$ is the kernel of $\pi_i\cmp f$, where $\pi_i$ is the canonical projection $\map {\pi_i}G{J_i\by{L_i}}$.
\end{pf}

The above lemma fails for disconnected groups.

\begin{ex}\label{woiejt}
Let $\T$ denote the circle group and let $\map \al{\T}{\T}$ be the 
conjugacy automorphism, i.e. $\al(z)=\ovr z$. Let $H$ be the $2$-element group generated by $\al$, i.e. $H=\dn1\al$, where $1$ denotes the identity.
Finally, let $G=\T\amalgam H$ be the semidirect product of $\T$ by $H$.
That is, $G=\T\times H$ and the group operation is
$$(z,h) \cdot (z',h') = (z h(z'), hh').$$
We claim that all proper normal subgroups of $G$ are contained in $\T\times\sn1$.

For fix $N\normalsub G$ and assume $(a,\al)\in N$ for some $a\in \T$. Fix $(b,\al)\in G$ and let $z\in \T$ be such that $az^2 = b$. Then $(z,1)\cdot (a,\al)\cdot (z^{-1},1)\in N$. On the other hand,
$$(z,1)\cdot (a,\al)\cdot (z,1)^{-1} = (za,\al)\cdot(z^{-1},1) = (za\ovr z^{-1},\al) = (z^2a,\al) = (b,\al).$$
It follows that $\T\times\sn\al\subs N$.
Further,
$(z,\al)\cdot(1,\al) = (z,1)$,
which shows that also $\T\times\sn 1\subs N$ and consequently $N=G$.

Thus, $G$ has no nontrivial homomorphisms onto infinite simple compact Lie groups.
Finally, supposing $G$ has a normal resolution
$$G=N_0\normalsup N_1 \normalsup \dots \normalsup N_n \normalsup N_{n+1}=\sn1$$
satisfying the assertion of Lemma~\ref{wefjwei},
we conclude that $N_i\by {N_{i+1}}$ must be finite for each $i\loe n$. Assuming $N_i\ne N_{i+1}$ for $i<n+1$, we see that $N_1=\T$ and $N_2$ must be a finite subgroup of $\T$, therefore $N_1\by N_2\iso \T$ is infinite. This is a contradiction.
\end{ex}

\begin{lm}\label{glowka}
Let $G$ be a connected compact group, let $K\normalsub A\normalsub G$. Then there exists $B\normalsub G$ such that $K\normalsub B\normalsub A\normalsub G$ and either $A\by B$ is finite or else there exists $C\normalsub G$ satisfying $A\cap C = B$ and $G\by C$ is a simple compact Lie group.
\end{lm}

\begin{pf}
Let $\En$ be the family of all $N\normalsub G$ such that $K\subs N$ and $G\by N$ is a Lie group. By Fact~\ref{f1} applied to $G\by K$, $\bigcap\En=K$.

Fix $N\in\En$ so that $N\not\sups A$.
Let $\sett{H_i}{i\loe n}$ be a normal resolution of $G\by N$ obtained from Lemma~\ref{wefjwei}. Let $N_i\normalsub G$ be such that $G\by{N_i} = (G\by N)\by {H_i}$.
Note that
$$G=N_0\normalsup N_1 \normalsup \dots \normalsup N_n=N.$$
Let $j\loe n$ be minimal such that $N_j\not\sups A$. Then $j>0$ and $N_{j-1}\sups A$. Let $B=A\cap N_j$. Clearly, $K\subs B$. We claim that $B$ is as required.

For suppose first that $H_{j-1}\by{H_j}$ is finite. Then $A\by B = (N_{j-1}\cap A)\by{(N_j\cap A)}$ is finite too.

Suppose now that there is $D\normalsub G\by N$ such that $(G\by N)\by D$ is a simple Lie group and $D\cap H_{j-1}=H_j$. 
Denote by $\map q{G}{G\by N}$ the canonical quotient map and let $C=\inv qD$. Then $G\by C \iso (G\by N)\by D$ is a simple Lie group. Finally, $C\cap A = B$. 
\end{pf}

Before stating the main result, we formulate and prove two simple lemmas about groups. The second one should clarify the meaning of Lemma~\ref{glowka}.

\begin{lm}\label{box1}
Let $N,L$ be two groups and let $G$ be a subgroup of $N\times L$. 
Further, let $q=\pr_N\rest G$, where $\pr_N, \pr_L$ denote the projections onto $N$ and $L$ respectively.
If $\img qG=N$ and $\img{\pr_L}{\ker(q)}=L$ then $G=N\times L$.
\end{lm}

\begin{pf}
The assumptions say, in particular, that $\sn 1 \times L \subs G$. Fix $(x,y)\in N\times L$. 
Since $q$ is an epimorphism, there is $z\in L$ such that $(x,z)\in G$. Now observe that $(1,z^{-1}y)\in G$. Hence $(x,y)=(x,z)\cdot(1,z^{-1}y)\in G$.
\end{pf}

\begin{lm}\label{giot}
Let $G$ be a compact group and let $A,B,C$ be closed normal subgroups of $G$ such that $A\cap C = B$ and $G\by C$ is a simple compact Lie group. If $A\by B$ is infinite then the canonical epimorphism $\map q{G\by B}{G\by A}$ is isomorphic, in the category of topological groups, to the projection $\map{\pr_{G\by A}}{(G\by A)\times(G\by C)}{G\by A}$.
\end{lm}

\begin{pf}
Let $\map p{G\by B}{G\by C}$ be the other canonical epimorphism (which is defined because $B\normalsub C$).
The equation $A\cap C=B$ means that the diagonal map $q\diag p$ (defined by $(q\diag p)(x)=(q(x),p(x))$) is a monomorphism (therefore also a topological embedding) of $G\by B$ into the product $(G\by A)\times(G\by C)$.
Thus, we may assume that $G\by B\subs (G\by A)\times(G\by C)$ and $q,p$ are the respective projections restricted to $G\by B$.
Now, $\ker(q)=A\by B$ is infinite and $p\rest(A\by B)$ is one-to-one, therefore $\img p {A\by B}$ is an infinite closed normal subgroup of $G\by C$.
Since $G\by C$ is a simple compact Lie group, we conclude that $\img p{\ker q}=G\by C$.
Finally, Lemma~\ref{box1} gives that $G\by B=(G\by A)\times(G\by C)$.
\end{pf}

\begin{tw}\label{ghorwgeoeqw}
Let $G$ be a connected compact group. Then $G=\liminv\S$, where $\S=\invsys Gp\vartheta \al\beta$ is a continuous transfinite inverse sequence in the category of compact groups such that $\vartheta$ is an ordinal, $G_0$ is the trivial group, and for each $\al<\vartheta$ either $p^{\al+1}_\al$ is an epimorphism with finite kernel or else $p^{\al+1}_\al$ is isomorphic to a projection of the form $\map \pr {G_\al\times H}G_\al$, where $H$ is a simple compact Lie group.

Furthermore, $\vartheta$ is a cardinal and
\begin{enumerate}
	\item[(a)] $\vartheta$ is finite if and only if $G$ is a Lie group.
	\item[(b)] $\vartheta=\omega$ if and only if $G$ is metrizable and not a Lie group.
	\item[(c)] $\vartheta=\w(G)$ as long as $G$ is not a Lie group.
\end{enumerate}
\end{tw}

\begin{pf}
Let $\sett{H_\al}{\al<\vartheta}$ be a normal resolution in $G$ such that for each $\al<\vartheta$
either $H_\al\by{H_{\al+1}}$ is finite
or else there exists $C_\al\normalsub G$ such that $C_\al\cap H_{\al}=H_{\al+1}$ and $G\by {C_\al}$ is a simple Lie group.
The existence of such a resolution follows, by transfinite induction, from Lemma~\ref{glowka}.
Indeed, having a partial normal resolution $\sett{H_\xi}{\xi<\al}$ like above, we set $H_\al=\bigcap_{\xi<\al}H_\xi$ in case $\al$ is a limit ordinal,
and we set $H_{\al}$ to be the group $B$ obtained from Lemma~\ref{glowka} with $A=H_{\al-1}$ and $K=\sn1$.
We finish this construction when $H_\al$ becomes trivial (and we set $\vartheta=\al$).

The normal resolution can be translated to a continuous inverse sequence in the category of topological groups, whose limit is $G$.
Namely, $G=\liminv\S$, where $\S=\invsys Gp\vartheta \al\beta$, 
$G_\al =G\by {H_\al}$ and $\map{p^\beta_\al}{G_\beta}{G_\al}$ is the canonical epimorphism (which is well defined because $H_\beta\normalsub H_\al$ for $\al<\beta$).
It remains to check that $\S$ has the desired properties. 

Note that $\S$ is continuous, because $H_\delta=\bigcap_{\xi<\delta}H_\xi$ for every limit ordinal $\delta\loe\vartheta$.

Fix $\al<\vartheta$. Observe that $\ker(p^{\al+1}_\al)=H_{\al}\by{H_{\al+1}}$.
Thus, if $H_{\al}\by{H_{\al+1}}$ is finite, $p^{\al+1}_\al$ is an epimorphism with finite kernel.
Otherwise,
there exists $C\normalsub G$ such that $G\by C$ is a simple compact Lie group and $C\cap H_\al = H_{\al+1}$.
By Lemma~\ref{giot}, $G\by{H_{\al+1}}$ is isomorphic to $(G\by{H_\al})\times (G\by C)$ and $p^{\al+1}_\al$ is isomorphic to the projection $\pr_{G\by H_\al}$ onto $G\by {H_\al}$.

This completes the proof of the first part. 
So far, we only know that $\vartheta$ is an ordinal.
{From} now on, we assume that the sequence $\S$ is nowhere trivial, that is, $\ker(p^{\al+1}_\al)$ is nontrivial for every $\al<\vartheta$.
 It is easy to see (by transfinite induction) that $\w(G)=|\vartheta|$ whenever $\w(G)$ is uncountable.
Thus (c) holds if $G$ is not metrizable.

If $\vartheta$ is finite then $G$ is a Lie group.
This follows from the fact that if $H$ is a compact Lie group, $\map qGH$ is a continuous epimorphism and
either $\ker q$ is finite
or $G$ is isomorphic to $H\times L$, where $L$ is a Lie group,
then $G$ is a compact Lie group too.

Now suppose $\vartheta$ is infinite.
We claim that $G$ is not a Lie group.
It suffices to show it in case where $\vartheta=\omega$, since otherwise $G_\vartheta$ is a quotient of $G$.
So assume $\vartheta=\omega$ and suppose $V$ is a neighborhood of the identity of $G$ that contains no nontrivial subgroups of $G$.
Let $\map{p_n}G{G_n}$ denote the projection.
Sets of the form $\inv{p_n}U$, where $\ntr$ and $U$ is open in $G_n$, form a basis for the topology of $G$,
therefore we can find $k\in\omega$ and a neighborhood $W$ of the identity of $G_k$ such that $\inv{p_k}W\subs V$.
Now, $K=\ker(p^{k+1}_k)$ is nontrivial and contained in $W'=\inv{(p^{k+1}_k)}W$.
Finally, $\inv{p_{k+1}}K$ is a nontrivial subgroup of $G$ contained in
$$\inv {p_{k+1}}{W'} = \inv{p_{k+1}}{\inv{(p^{k+1}_k)}W} = \inv {p_k}W \subs V,$$
a contradiction.

This shows (a).
In order to show (b) and to complete the proof of (c), it remains to show that $\vartheta\loe\omega$  whenever $G$ is metrizable. Formally, from the above proof we can conclude that $\vartheta$ is a countable infinite ordinal whenever $G$ is a metrizable group which is not a Lie group. We need to modify the arguments in this case to obtain that $\vartheta=\omega$.

Namely, assume that $\w(G)=\aleph_0$ and that $G$ is not a Lie group.
Start by embedding $G$ into $\prod_{\ntr}L_n$, where each $L_n$ is a compact connected Lie group.
We also assume that each projection restricted to $G$ is an epimorphism.
This induces a normal resolution
$$G = N_0\normalsup N_1\normalsup \dots\normalsup \sn1$$
of length $\omega$, such that for each $\ntr$  there exists $C_n\normalsub G$ with $G\by C_n$ a Lie group and $C_n\cap N_n=N_{n+1}$.

Fix $k\in\nat$. Using Lemma~\ref{glowka} and (possibly transfinite) induction, we obtain a strictly decreasing sequence $\sett{H_i}{i<\zeta}$ of normal subgroups of $G$ such that $H_0=N_k$, $\bigcap_{i<\zeta}H_i=N_{k+1}$ and each pair $H_i,H_{i+1}$ satisfies the assertion of Lemma~\ref{glowka}.
That is, either $H_i\by{H_{i+1}}$ is finite or there exists $C\normalsub G$ such that $G\by C$ is a simple compact Lie group and $C\cap H_i=H_{i+1}$.

We claim that the sequence $\sett{H_i}{i<\zeta}$ is finite (i.e. $\zeta$ is a natural number).
After we show this, we can modify the original resolution $\sett{N_n}{\ntr}$ by `inserting' a finite sequence between each $N_n,N_{n+1}$
so that the resulting resolution (still of length $\omega$) will induce an inverse sequence with all required properties, like it was done before.

Fix $C\normalsub G$ so that $G\by C$ is a Lie group and $C\cap N_k=N_{k+1}$. Let $\map qG{G\by C}$ be the canonical quotient map. Then $\img q{H_i}$ is a normal subgroup of $G\by C$. Thus we have a normal sequence
$$G\by C \normalsup \img q{H_0} \normalsup \img q{H_1} \normalsup \dots,$$
in the Lie group $G\by C$.
Observe that $\img q{H_i}\ne \img q{H_{i+1}}$ for each $i$.
Suppose otherwise and choose $x\in H_i\ne H_{i+1}$ and $y\in H_{i+1}$ so that $xC = yC$.
This means that $x = yc$ for some $c\in C$.
On the other hand, $xy^{-1}\in N_k$, so $c\in N_k\cap C=N_{k+1}$.
Thus, $x\in H_{i+1}N_{k+1}\subs H_{i+1}$, because $N_{k+1}\subs H_{i+1}$, a contradiction.
By this way we have constructed a strictly decreasing sequence $\sett{\img q{H_i}}{i<\zeta}$ of closed normal subgroups of a compact Lie group.
By arguments from the proof of (b), we conclude that $\zeta$ is finite.

This completes the proof of the theorem.
\end{pf}

\begin{uwga}
Example~\ref{woiejt} witnesses that Theorem~\ref{ghorwgeoeqw} fails for disconnected compact Lie groups.
However, Theorem~\ref{ghorwgeoeqw} is true both for $0$-dimensional and for Abelian groups.

In fact, if $G$ is a $0$-dimensional compact group then $G$ embeds into $\prod_{i<\kappa}F_i$, where each $F_i$ is a finite group, see~\cite[Theorem~1.34]{HM}.
This embedding gives rise to a continuous inverse sequence $\S=\invsys Gp\kappa\al\beta$ with limit $G$,
where $G_\al$ is the projection of $G$ onto $\prod_{i<\al}F_i$
and $p^{\al+1}_\al$ is the restriction of the projection $\map \pr{G_\al\times F_\al}{G_\al}$. 
Clearly, $\ker(p^{\al+1}_\al)$ is finite, because it embeds into $F_\al$.

Now, if $G$ is an Abelian compact group, then $G$ is topologically isomorphic to a subgroup of $\T^\kappa$, where $\T$ denotes the circle group, see e.g.~\cite[Corollary 2.31]{HM}.
Again, this induces a continuous inverse sequence $\S=\invsys Gp\kappa\al\beta$ with limit $G$,
where $G_\al$ is the projection of $G$ onto $\T^\al$
and $p^{\al+1}_\al$ is the restriction of the projection $\map \pr{G_\al\times \T}{G_\al}$.
Since $\T$ is a simple compact Lie group, if $\ker(p^{\al+1}_\al)$ is infinite then $G_{\al+1}=G_\al\times \T$ and $p^{\al+1}_\al$ is the standard projection.
\end{uwga}

By the above remark, we get

\begin{prop}\label{erhoou}
Let $G$ be a $0$-dimensional compact group.
Then $G=\liminv\S$, where $\S=\invsys Gp\kappa\al\beta$ is a continuous inverse sequence in the category of topological groups, $G_0$ is the trivial group and each $p^{\al+1}_\al$ is an epimorphism with a finite kernel.

Furthermore, if $G$ is infinite then $\kappa=\w(G)$.
\end{prop}

\begin{prop}\label{abcpaie}
Let $G$ be an Abelian compact group.
Then $G=\liminv\S$, where $\S=\invsys Gp\kappa\al\beta$ is a continuous inverse sequence in the category of topological groups, $G_0$ is the trivial group and for each $\al<\kappa$, either $p^{\al+1}_\al$ is an epimorphism with a finite kernel or, up to isomorphism, $G_{\al+1}=G_\al\times\T$ and $p^{\al+1}_\al$ is the standard projection.

Furthermore, if $G$ is not a Lie group then $\kappa=\w(G)$.
\end{prop}

\begin{uwga}
Theorem~\ref{ghorwgeoeqw}, as well as the above propositions, can be generalized in the following way:

\begin{claim}\label{klejm2341}
Given a continuous epimorphism $\map fG{G'}$ of connected compact groups, there exists a continuous inverse sequence $\S=\invsys Gp \vartheta \al\beta$ satisfying all the assertions of Theorem~\ref{ghorwgeoeqw} (in particular $G=\liminv\S$) except that $G_0=G'$ and the projection $\map {p_0}G{G_0}$ is isomorphic to $f$.
\end{claim}

Similar statement holds both for Abelian and for $0$-dimensional compact groups.

Claim~\ref{klejm2341} can be proved by slight modifications of the arguments in the proof of Theorem~\ref{ghorwgeoeqw} or by using the following observation:
Assume $\sett{H_i}{i<\vartheta}$ is a normal resolution in $G'$ satisfying the assertions 
of Lemma~\ref{wefjwei} (except that $\vartheta$ may be infinite) and let $L_i=\inv f{H_i}$, $i<\vartheta$.
Then $\sett{L_i}{i<\vartheta}$ is a continuous normal sequence in $G$ that also satisfies the assertions of Lemma~\ref{wefjwei}, 
that is, for each $i<\vartheta$ either $L_i \by L_{i+1}$ is finite or else there exists $D\normalsub G$ such that $D\cap L_i = L_{i+1}$ and $G\by D$ is a simple compact Lie group.
Note that $\sett{L_i}{i<\vartheta}$ induces the inverse sequence asked by Claim~\ref{klejm2341}.
\end{uwga}

\section{Supercompactness}

Recall from the introduction that a closed subbase $\Bee$ of a topological space $X$ is \emph{binary} if $\bigcap \Aaa\nnempty$ whenever $\Aaa\subs \Bee$ is linked, i.e. $A_0\cap A_1\nnempty$ for every $A_0,A_1\in \Aaa$.
It is easy to see that, given a binary subbase $\Bee$, the family
$$\Bee^\cap=\setof{\bigcap\Ef}{\Ef\subs \Bee}$$
is binary too.

The following two statements belong to the folklore.

\begin{lm}\label{zapwoefhoew}
Let $\En$ be a family of closed subsets of a compact space $X$ and assume $\En$ is stable under finite intersections. Then $\En$ is a closed subbase if and only if for every $x\in X$ and for every open $U\subs X$ with $x\in U$, there is a finite $\Aaa\subs\En$ satisfying
$$x\in\Int(\bigcup\Aaa)\subs \bigcup\Aaa \subs U.$$
\end{lm}

\begin{pf}
Suppose first that $\En$ is a closed subbase. Fix $x\in U=\Int U$ and choose an open set $V$ so that $x\in V\subs \cl V\subs U$. Then $\cl V$ is the intersection of some family $\Cee$ that consists of finite unions of members of $\En$. By compactness, there is a finite $\Cee_0\subs \Cee$ such that $\bigcap\Cee_0\subs U$. 
Using the law of distributivity for sets, we see that $\bigcap\Cee_0=A_0\cup\dots\cup A_{n-1}$, where each $A_i$ is a finite intersection of members of $\En$. By our assumption, $A_i\in\En$ for $i<n$.

Suppose now that $\En$ satisfies the assertion of the lemma. We need to show that the family $\En'$ consisting of all finite unions of elements of $\En$ is a base for closed sets. Fix $x\in X$ and a closed set $B\subs X\setminus \sn x$. For each $y\in B$ choose a finite $\En_y\subs \En$ such that $y\in\Int(\bigcup\En_y)$ and $x\notin \bigcup\En_y$. Let $S\subs B$ be finite such that $\setof{\bigcup\En_y}{y\in S}$ covers $B$ (here we use compactness again). Then $\Ef=\bigcup_{y\in S}\En_y\in\En'$ and $x\notin \bigcup\Ef\sups B$.
\end{pf}

\begin{lm}\label{cov}
Assume $\Bee$ is a binary subbase in a compact space $X$, stable under finite intersections. Let $\Yu$ be an open cover of $X$. Then
$$\Bee'=\setof{B\in\Bee}{(\exists\;U\in\Yu)\; B\subs U}$$
is a binary subbase for $X$.
\end{lm}

\begin{pf}
It is obvious that $\Bee'$ is binary. The fact that it is a closed subbase follows immediately from Lemma~\ref{zapwoefhoew}.
\end{pf}

The next result should be known, although it was probably never stated explicitely in the literature. Its special case for compact groups was actually used by Mills in~\cite{Mills}, 

\begin{tw}\label{etwoqsd}
Let $K=\liminv \S$ where $\S=\invsys Kp\kappa \al\beta$ is a continuous inverse sequence of compact spaces with surjections such that the following conditions are satisfied.
\begin{enumerate}
	\item[(i)] $K_0$ is supercompact.
	\item[(ii)] For every $\al<\kappa$, either $\map{p^{\al+1}_{\al}}{K_{\al+1}}{K_\al}$ is a local homeomorphism or else it is homeomorphic to a projection $\map{\pr}{K_\al\times X}{K_\al}$, where $X$ is a supercompact space.
\end{enumerate}
Then $K$ is a supercompact space.
\end{tw}

\begin{pf}
We shall prove by induction the following statement: given a closed binary subbase $\Aaa_0$ for $K_0$, stable under finite intersections, there exists a closed binary subbase $\Aaa_\kappa$ for $K$, also stable under finite intersections and such that $\inv{p_0}A\in \Aaa_\kappa$ for every $A\in \Aaa_0$.
We shall say that $\Aaa_\kappa$ is an \emph{extension} of $\Aaa_0$ with respect to the surjection $p_0$. Recall that $p_\al$ denotes the projection from $K$ onto $K_\al$ ($\al<\kappa$).

We construct inductively families $\Aaa_\al$ such that $\Aaa_\al$ is a closed binary subbase for $K_\al$, stable under finite intersections, and $\Aaa_\beta$ is an extension of $\Aaa_\al$ with respect to $p^{\beta}_\al$, whenever $\al<\beta$.

Fix an ordinal $\delta\loe \kappa$ and suppose the families $\Aaa_\al$ have been constructed for $\al<\delta$. If $\delta$ is a limit ordinal, we set
$$\Aaa_\delta=\setof{\inv{(p^{\delta}_\al)}A}{A\in\Aaa_\al,\; \al<\delta},$$
where $p^\kappa_\al=p_\al$ (in case $\delta=\kappa$).
It is straight that $\Aaa_\delta$ is a closed binary subbase for $K_\delta$ (or for $K$, in case $\delta=\kappa$).

Suppose now that $\delta = \al+1$ and let $f=p^{\al+1}_\al$.
We consider two cases.

\paragraph{Case 1} $f$ is a local homeomorphism.

Let $\Yu$ be an open cover of $K_{\al+1}$ such that $f\rest U$ is a homeomorphism and $\img fU$ is open in $K_\al$ for every $U\in\Yu$  . Let $\Vee$ be a star-refinement of $\Yu$ and define
$$\Aaa_{\al+1} = \setof{\inv fA}{A\in \Aaa_\al}\cup\setof{B\subs K_{\al+1}}{\img fB\in\Aaa_\al \land (\exists\;V\in\Vee)\; B\subs V}.$$
The equation $\img f{\inv fA \cap B} = A \cap  \img fB$ implies that $\Aaa_{\al+1}$ is stable under finite intersections.
By Lemma~\ref{zapwoefhoew} and Lemma~\ref{cov} it is a closed subbase for $K_{\al+1}$.
It is remains to check that $\Aaa_{\al+1}$ is binary.
Clearly, $\setof{\inv fA}{A\in \Aaa_\al}$ is binary.
Given a linked family $\El$ which refines $\Vee$ and 
such that $\img fB\in\Aaa_\al$ for each $B\in\El$, notice that $\bigcup\El$ is contained in some $U\in\Yu$, because $\Vee$ is a star-refinement of $\Yu$.
It follows that $f$ is one-to-one on $\bigcup\El$ and consequently $\bigcap\El\nnempty$, because $\bigcap_{B\in\El}\img fB\nnempty$.
Now, the proof that $\Aaa_{\al+1}$ is binary can be reduced to the case of a linked family of the form
$$\{\inv fA, B_0,B_1\},$$
where $A\in\Aaa_\al$, $\img f{B_i}\in\Aaa_\al$ and  $B_i\subs V_i$ for some $V_i\in\Vee$, $i<2$. Since $B_0\cap B_1\nnempty$, again using the fact that $\Vee$ is a star-refinement of $\Yu$, there exists $U\in\Yu$ such that $B_0\cup B_1\subs V_0\cup V_1 \subs U$. Let $C=\inv fA \cap B_0$. Then $\img fC=A\cap \img f{B_0}\in\Aaa_\al$ and $C\subs V_0\subs U$. Now, if $C\cap B_1=\emptyset$ then $\img fC\cap \img f{B_1}=\emptyset$, because $f\rest U$ is one-to-one. But this contradicts the fact that the family $\{A,\img f{B_0}, \img f{B_1}\}\subs \Aaa_\al$ is linked.

\paragraph{Case 2} $f$ is homeomorphic to a projection $\map \pi{K_{\al}\times X}{K_\al}$, where $X$ is supercompact.

In particular, up to homeomorphism, $K_{\al+1}=K_\al\times X$.
Let $\Bee$ be a closed binary subbase for $X$, stable under finite intersections. Define
$$\Aaa_{\al+1} = \setof{A\times B}{A\in\Aaa_\al,\; B\in\Bee}.$$
It is clear that $\Aaa_{\al+1}$ is a closed binary subbase for $K_{\al+1}$, stable under finite intersections. Assuming $X\in\Bee$, we are sure that $\Aaa_{\al+1}$ is an extension of $\Aaa_\al$ with respect to $f$.

This completes the proof.
\end{pf}

\begin{wn}[Mills]
Every compact group is supercompact.
\end{wn}

\begin{pf}
Let $G$ be a compact group. Then $G$ is homeomorphic to $G_0\times Z$, where $G_0$ is a connected compact group and $Z$ is a $0$-dimensional compact group, necessarily homeomorphic to a Cantor cube (see \cite {HM} for details). Since Cantor cubes are obviously supercompact and any product of supercompact spaces is supercompact, we may assume that $G$ is connected.

By Theorem~\ref{ghorwgeoeqw}, $G$ is the limit of a continuous inverse sequence satisfying the conditions of Theorem~\ref{etwoqsd}. By the latter theorem, $G$ is supercompact.
\end{pf}

A family of $\Ef$ in a group $G$ is \emph{translation invariant} if $gF\in\Ef$ whenever $g\in G$ and $F\in\Ef$.
Note that the circle group $\T$ has a translation invariant closed binary subbase, consisting of all closed arcs of length $<2\pi/3$. Obviously, every finite group has a translation invariant binary subbase consisting of singletons.
Thus, using Propositions~\ref{abcpaie}, \ref{erhoou} and examining the proof of Theorem~\ref{etwoqsd}, we obtain

\begin{wn}
All compact Abelian groups as well as all $0$-dimensional compact groups have translation invariant closed binary subbases.
\end{wn}

We do not know whether every simple compact Lie group has a translation invariant closed binary subbase. If so, then the above corollary would be valid for all connected compact groups.

\begin{uwga}
We make some comments about Mills' work.
First of all, his proof of supercompactness of compact groups also heavily used Facts~\ref{f1} and \ref{f2} and next that supercompactness is preserved by local homeomorphisms.
Our Theorem~\ref{ghorwgeoeqw} simplifies one of the crucial steps in Mills' proof, where originally some lifting was needed. 

Actually, the main result of Mills~\cite{Mills} is the following extension theorem.

\begin{tw}
Let $\map fG{G'}$ be a continuous epimorphism of compact groups and let $\Aaa$ be a binary closed subbase for $G'$ . Then there exists a closed binary subbase $\Bee$ for $G$ such that $\inv fA\in\Bee$ for every $A\in\Aaa$.
\end{tw}

Taking $G'=\sn1$, one gets the conclusion that every compact group is supercompact.

In case $G$ is connected (or Abelian, or $0$-dimensional), the above result follows easily from Theorem~\ref{etwoqsd} together with Claim~\ref{klejm2341} (or its version corresponding to Abelian/$0$-dimensional groups).

For a general $G$, the argument of Mills can be sketched as follows.

First, using Fact~\ref{f1}, find a continuous inverse sequence $\S=\invsys Gp\kappa\al\beta$ in the category of topological groups, such that $G=\liminv\S$, $G_0=G'$, $p_0=f$ and for each $\al<\kappa$ there exists a compact Lie group $H_\al$ such that,
up to isomorphism, $G_{\al+1}$ is a subgroup of $G_\al\times H_\al$ and $p^{\al+1}_\al$ is the restriction to $G_{\al+1}$ of the projection $\map{\pr_{G_\al}}{G_\al\times H_\al}{G_\al}$.
We may also assume that $\img {\pr_{H_\al}}{G_{\al+1}} = H_\al$. 
We extend a given closed binary subbase $\Aaa_0$ in $G_0=G'$ to a binary subbase
$\Aaa_\al$ in $G_\al$ for each $\al<\kappa$ by induction on $\al$.
The limit step makes no problems, so it suffices to consider the successor step.
Thus, suppose $\Aaa_\al$ has already been defined.

If $H_\al$ is connected, actually the arguments from the proof of Theorem~\ref{ghorwgeoeqw} still apply: we can ``insert" a finite sequence between $G_\al$ and $G_{\al+1}$ satisfying the assertions of Theorem~\ref{ghorwgeoeqw}. Further, applying the proof of Theorem~\ref{etwoqsd}, we see that the closed binary subbase $\Aaa_\al$ for $G_\al$ can be extended to a closed binary subbase $\Aaa_{\al+1}$ for $G_{\al+1}$.

In case $H_\al$ is disconnected, we use the above argument for its identity component $H_\al^0$ and then we apply the fact that, as a topological space, $H_\al=H_\al^0\times F$, where $F$ is a finite set.

We leave the verification of the details to interested readers.
\end{uwga}

\end{document}